\documentclass[10pt,twoside]{article}
\usepackage{graphicx}
\usepackage{amsmath}
\usepackage{Latex-document}
\usepackage{amssymb}

\title{\bf  Finite Dimensional Approximations \vskip -2mm in Geometry\vskip 6mm}

\author{M. Furuta\vspace*{-0.5cm}\thanks{Graduate School of Mathematical Sciences, University of Tokyo,
3-8-1, Komaba, Meguro-ku, Tokyo 153-8914, Japan. E-mail: furuta@ms.u-tokyo.ac.jp}}

\date{\vspace{-8mm}}

\newcommand{\HH}{{\mathbb{H}}}

\newcommand{\Z}{{\mathbb Z}}

\newcommand{\R}{{\mathbb{R}}}

\newcommand{\C}{{\mathbb{C}}}
\newcommand{\Q}{{\mathbb{Q}}}

\newcommand{\Ker}{{\mathrm Ker}}
\newcommand{\Coker}{{\mathrm Coker}}

\newcommand{\ind}{{\mathrm ind}\,}
\newcommand{\sign}{{\mathrm sign}\,}

\begin{document}

\maketitle

\thispagestyle{first} \setcounter{page}{395}

\begin{abstract}\vskip 3mm
In low dimensional topology,
we have some invariants defined by using
solutions of some nonlinear elliptic operators.
The invariants could be understood
as Euler class or  degree in
the ordinary cohomology,
in infinite dimensional setting.
Instead of looking at the solutions,
if we can regard  some kind of
homotopy class of the operator itself
as an invariant, then the refined version
of the invariant is understood
as Euler class or degree in
cohomotopy theory.
This idea can be carried out for
the Seiberg-Witten equation
on $4$-dimensional manifolds
and we have some applications
to $4$-dimensional topology.

\vskip 4.5mm


\end{abstract}

\vskip 12mm

\section{Introduction} \label{section 1}\setzero
\vskip-5mm \hspace{5mm}

The purpose of this paper is to review
the recent developments in
a formal framework to extract topological
information from nonlinear elliptic operators.

We also explain some applications of the idea
to 4-dimensional topology by using
the Seiberg-Witten theory.


A prototype is
the notion of index for linear elliptic operators.
In this introduction we explain this linear case.
Later we mainly explain the Seiberg-Witten case.

Let $D:\Gamma(E^0) \to \Gamma(E^1)$ be an
linear elliptic operator  on a close
manifold $X$. The index $\ind D$ is defined to be
$$
\ind D =\dim \Ker D - \dim \Coker D.
$$
We can extend this definition as follows.
Take any decomposition
$D=L\oplus L': V^0 \oplus W^0 \to V^1 \oplus W^1$.
such that $L:V^0 \to V^1$ is a linear map
between two finite dimensional
vector spaces and that $L':W^0 \to W^1$ is
an isomorphism between infinite dimensional vector
spaces.
Then we have
$$
\ind D = \dim V^0 - \dim V^1.
$$
It is easy to check that the right-hand-side
is independent of the choice of the decomposition.
For example we have decomposition satisfying
$V^0=\Ker D$, $V^1\cong \Coker D$, $L=0$, which
gives the former definition of the index.

An important property of $\ind D$ is its invariance
under continuous variation of $D$.
This property is closely related to the above
well-definedness.

Another way to understand this property is
to consider the whole space of Fredholm maps.
Then the given map $D$ sits in the space and
the $\ind D$ is nothing but the label of
the connected component containing $D$.

In other words, there are presumably
three possible attitudes:
\begin{enumerate}
\item
The essential data is ``supported" on $\Ker D$ and $\Coker D$.
\item
It is convenient to look at ``some" finite dimensional
approximation. $L:V^0 \to V^1$.
\item
The essential data is
the whole map $D:\Gamma(E^0)\to \Gamma(E^1)$.
\end{enumerate}

When one considers a family of elliptic operators and
tries to define the index of the family,
it is not enough to look at their kernels and cokernels.

It is tempting to regard the finite dimensional
approximation as a topological version of the notion
of ``low energy effective theory" in physics.
In this story, the whole map $D$ would be
regarded as a given original theory.

In this paper we explain a nonlinear version
of the notion of index which is formulated
by using finite dimensional approximations.

\section{Non-linear cases} \label{section 2}
\setzero\vskip-5mm \hspace{5mm}

While every elliptic operator on a closed manifold
has its index as topological invariant,
it is quite rare that a nonlinear elliptic
operator gives some topological invariant.

We have three examples of this type of invariants:
the Donaldson invariant,
the Gromov-Ruan-Witten invariant and
the Seiberg-Witten invariant. Moreover,
the Casson invariant is regarded a variant of
the Donaldson invariant. Some
other finite type invariants for $3$-manifolds
are also supposed to be related to these kinds
of invariant \cite{SeibergWitten}.

Even for these cases, however,
it is not obvious how to proceed to obtain
nonlinear version of index in full generality.

Let us first give several examples of
finite dimensional approximations.

\begin{enumerate}
\item
C. Conley and E. Zehnder solved
the Arnold conjecture for torus by
reducing a certain variational problem
to a finite dimensional Morse theory \cite{ConleyZehnder}.
\item
Casson's definition of the Casson invariant.
Taubes gave an interpretation of
the Casson invariant via gauge theory \cite{Taubes}.
In other words, Casson's construction
gave  a finite dimensional approximation
of the gauge theoretical setting by Taubes.
(The statement of the Atiyah-Floer conjecture
could be regarded as a {\it partial} finite
dimensional approximation along fibers.)

\item Seiberg-Witten equation.
The moduli space of Seiberg-Witten equation
is known to be compact for closed $4$-manifolds.
This enables us to globalize the Kuranishi
construction to obtain finite dimensional
approximations \cite{Furuta}, \cite{BauerFuruta}.
\item Seiberg-Witten-Floer theory
C. Manolescu and P.B. Kronheimer
defined Floer homotopy type for
Seiberg-Witten theory,
which is formulated as spectrum
\cite{Manolescu} \cite{KronheimerManolescu}.
\item
Kontsevish explained
an idea to define invariants of $3$-manifolds
by using configuration spaces.
This idea was realized by
Fukaya \cite{Fukaya},
Bott-Cattaneo \cite{BottCattaneo1} \cite{BottCattaneo2},
and Kuperberg-Thurston \cite{KuperbergThurston}.
Formally the configuration spaces
appear as finite approximations
of certain path spaces.
\end{enumerate}

\section{Kuranishi construction} \label{section 3} \setzero\vskip-5mm \hspace{5mm}

While the index is regarded as the infinitesimal
information of a nonlinear elliptic operator,
its local information is given by the Kuranishi
map, which has been used to describe local structure
in various moduli problems \cite{Kuranishi}.

A few years ago
the Arnold conjecture was
solved in a fairly general setting
and the Gromov-Ruan-Witten invariant
was defined for general symplectic manifolds.
These works were done by
 several groups independently \cite{FukayaOno},
\cite{LiuTian}, \cite{Ruan1}, \cite{Siebert}.
A key of their arguments was
to construct virtual moduli cycle over $\Q$.

In their case, the point is to glue local structure
to obtain some global data to define
invariants. Since their invariants are
defined by evaluating cohomology classes,
it was enough to have the virtual moduli cycle.






\section{Global approximation} \label{section 4} \setzero\vskip-5mm \hspace{5mm}

The notion of Fukaya-Ono's Kuranishi structure or
Ruan's virtual neighborhood
is defined as equivalence class of
collections of maps, which define
the moduli space.
The collection of maps is necessary because
the moduli space as topological space
is not enough to recover
the nature of the singularity on it.

The data depend on the choice of various
choice of auxiliary data. When
we change the data, the change of
the moduli space is supposed to be given by a cobordism,
even with the extra structure we have to look at.

Suppose we would like to regard this structure itself
as an invariant. Then we have to
identify the place where the invariant lives.
Since cobordism classes are identified by
Pontrjagin-Thom construction,
what we need would be a certain stable version of
Pontrjagin-Thom construction.

In the case of symplectic geometry
or Donaldson's theory, this construction
has not been done.
A main problem seems to describe
a finite dimensional approximation of the ambient space
where the compactification of the moduli space lies.
(The same problem occurs
for Kotschick-Morgan conjecture.)
Since the compactification is fairly complicated,
it is not straightforward to identify
the finite dimensional approximation.

However in the Seiberg-Witten case,
the moduli spaces are known to be compact
for closed $4$-manifolds and
it is not necessary to take any further
compactifications.

Let us briefly recall the Seiberg-Witten
equation for a closed Spin${}^c$ manifold $X$.
For simplicity we assume $b_1(X)=0$.
Let $W=W^0 \oplus  W^1$ be the spinor bundle and
${\cal A}$ be the space of connections on
$\det W^0 \cong \det W^1$.
Then the Seiberg-Witten equation is given by a map
$$
 \Gamma(W^0) \times {\cal A} \to
\Gamma(W^1) \times \Gamma(\wedge^+),
$$
where $\Gamma(\wedge^+)$ is the self-dual $2$-forms
for a fixed Riemannian metric.
This is an $U(1)$-equivariant map.
The inverse image of $0$ divided by $S^1$ is
the moduli space, which is known to be compact.

A finite approximation of the above map is
defined by
global version of the Kuranishi construction.
The approximation is a {\it proper}  $U(1)$-equivariant
map
$$
\C^{a_0} \oplus \R^{d_0} \to \C^{a_1} \oplus \R^{d_1}
$$
for some natural numbers $c_0,c_1,d_0$ and $d_1$.
The differences $c_0-c_1$ and $d_0-d_1$ depends only
on the topology of $X$ and its spin${}^c$-structure.

The invariant we have is
the stable homotopy class of the above $U(1)$-equivariant
proper map, or equivalently,
the $U(1)$-equivariant map from the sphere
$S(\C^{a_0} \oplus \R^{d_0})$ to the sphere
$S(\C^{a_1} \oplus \R^{d_1})$.


S. Bauer and the author pointed out that
the invariant constructed above
is a refinement of the usual
Seiberg-Witten invariant \cite{BauerFuruta}.

\section{4-dimensional topology and Seiberg-Witten \\ theory} \label{section 5} \setzero\vskip-5mm \hspace{5mm}

We explain some applications of
the finite dimensional approximation
to $4$-dimensional topology.

(1) Bauer's connected sum formula \cite{Bauer}

Suppose $X$ is the connected sum of $X_0$ and $X_1$. If the neck
of the connected sum is long enough, it is known that the moduli
space of the solution of the Seiberg-Witten equation (or
anti-self-dual equation) for $X$ is identified with the product of
the moduli spaces for $X_0$ and $X_1$. When $X_1=\overline{CP^2}$,
then this gives the blowing-up formula. When
$b^+(X_0),b^-(X_1)\geq 1$, this gives vanishing of the
Seiberg-Witten (or the Donaldson) invariant of $X=X_0 \# X_1$.
Bauer essentially showed that the product formula holds true for
the virtual neighborhood of the moduli spaces, if we use Ruan's
terminology. In the language of stable maps between spheres,
``product" becomes ``join". In particular Bauer's formula gives
the blowing-up formula for the refined invariant. When $b^+(X_0),
b^-(X_1)\geq 1$, the join is torsion. It is, however, not
necessary zero. In this way Bauer gave many new examples of
$4$-manifolds which are homeomorphic  but not diffeomorphic to
each other.

Ishida-Lebrun \cite{IshidaLebrun1} \cite{IshidaLebrun2}
obtained some applications of the
connected sum formula to Riemannian geometry.

(2) Intersection form of spin $4$-manifolds

When $4$-manifold is spin, we have
certain extra symmetry, and
the place where the invariant lives
is a set of $Pin(2)$-equivariant stable maps
\cite{Furuta}.


When $X$ is a closed spin $4$-manifold with $b_1(X)=0$,
the Seiberg-Witten map for the spin structure is
a $Pin(2)$-equivariant map formally given by
$$
\HH^\infty \oplus \tilde{\R}^\infty
\to\HH^\infty \oplus \tilde{\R}^\infty,
$$
where $\tilde{\R}$ is the  non-trivial
$1$-dimensional real representation space of $Pin(2)$.
and $\HH$ is the $4$-dimentional real
irreducible representation space of $Pin(2)$.
Let $\Z/4$ be the subgroup of $Pin(2)$ generated
by an element in $Pin(2)\setminus U(1)$.
The differences of the power $\infty$'s  are
given by the index of some elliptic operators.

A finite dimensional approximation is given by
a $Pin(2)$-equivariant proper map
$$
\HH^{c_0} \oplus \tilde{\R}^{d_0}
\to\HH^{c_1} \oplus \tilde{\R}^{d_1},
$$
for some $c_0,c_1,d_0,d_1$ satisfying
$$
c_0-c_1=-\frac{\sign(X)}{16},\,\,
d_0-d_1=b^+(X).
$$
This existence implies some inequality
between the signature and the second Betti number.

To obtain the inequality explicitly we can use
the following results.

{\bf Theorem}  \it
Suppose $k>0$ and $k \equiv a  \bmod 4$ for
$a=0,1,2,$ or $3$. Then
there does not exist a $G$-equivariant
continuous map from $S(\HH^{k+x}\oplus \tilde{\R}^{y})$
to
$S(\HH^x \oplus \tilde{\R}^{2k+a'-1+y})$.
for the following $G$ and $a'$.
\begin{enumerate}
\item
{\rm (B. Schmidt \cite{Schmidt} see
also \cite{Stolz} \cite{Crabb} \cite{Minami})}
$G=\Z/4$ and $a'=a$ for $a=1,2,3$.
\item
{\rm (F - Y.Kametani \cite{FurutaKametani2})}
$G=Pin(2)$ and $a'=3$ for $a=0$.
\end{enumerate}
\rm

From the above non-existence results, we
have the following inequality,
which is a partial result towards
the 11/8-conjecture $b^+ \geq 3|\sign(X)/16|$.

{\bf Theorem} \it
Let $X$ be a closed spin $4$-manifold with
$\sign(X)=-16k<0$. If $k \equiv a \bmod 4$
for $a=0,1,2$ or $3$, then we have
$ b^+ \geq 2k+b$, where $a'=a$ if $a=1,2,3$ and
$a'=3$ if $a=0$.
\rm

Equivariant version and $V$-manifold version
can be formulated similarly
\cite{Bryan}, \cite{Fang}, \cite{FukumotoFuruta},
\cite{Acosta}.
There are some applications of these extended versions:
\begin{enumerate}
\item
C. Bohr \cite{Bohr} and
R.~Lee - T.-J.~Li  \cite{LeeLi}
investigated the intersection forms of
closed even $4$-manifolds which
are not spin.
\item
Y.~Fukumoto, M.~Ue and the author
\cite{FukumotoFuruta} \cite{Fukumoto}
\cite{FukumotoFurutaUe} \cite{Ue},
and N.~Saveliev \cite{Saveliev}
investigated homology cobordims groups of
homology $3$-spheres.
\end{enumerate}
When $b_1>0$, we can construct
another closed spin $4$-manifold with $b_1=0$
without changing the intersection form.
It implies that we can assume $b_1=0$ to
obtain restriction on the intersection form.
However when the intersection form on $H^1(X)$
is non-trivial,
we may have a stronger restriction.
Y.~Kametani, H.~Matsue, N.~Minami and the author
found that such a phenomenon actually occurs if
there are $\alpha_1,\alpha_2,\alpha_3,\alpha_4
\in H^1(X,\Z)$ such that
$\langle \prod \alpha_i,[X]\rangle$
is odd \cite{FKMM}.

\section{Seiberg-Witten-Floer homotopy type}
\label{section 6} \setzero\vskip-5mm \hspace{5mm}

Recently
C.~Manolescu and P.~B.~Kronheimer extends
the above formulation for closed $4$-manifolds
to the relative version \cite{Manolescu},
\cite{KronheimerManolescu}.
Let us explain their theory briefly.

We mentioned that Conley-Zehnder
used a finite dimensional approximation
of a Morse function on an infinite dimensional space
to approach the Arnold conjecture for
torus. Following this line,
Conley exteded the notion
of Morse index and defined
the Conley index for compact isolated set \cite{Conley}.
The Conley index is not a number, but a
homotopy type of spaces.
Floer extracted some information from
the Conley index just by looking at
some finite dimensional skeleton of
the Conley index under some assumption.
Floer's formulation has the advantage that
the Floer homology is defined
even when the Conley index is not rigorously defined.

On the other hand
R.~L.~Cohen, J.~D.~S.~Jones and G.~B.~ Segal
tried to define certain stable homotopy type
directly which should be an extended version of
the Conley index \cite{CohenJonesSegal}.
They called it the Floer homotopy type.
At that time the Floer homology was defined only
for the Donaldson theory and
the Gromov-Ruan-Witten theory.
In these theories the moduli spaces are non-compact
in general.
This cause a serious difficulty to carry out their program.

In the Seiberg-Witten theory,
we have a strong compactness for the moduli spaces.
Manolescu and Kronheimer succeeded
to construct the Floer homotopy type
as spectra for the Seiberg-Witten theory by using
this compactness.

They also defined
relative invariant for $4$-manifolds with
boundary is also defined and it extends
the invariant in \cite{BauerFuruta}.

\section{Concluding remarks} \label{section 7} \setzero\vskip-5mm \hspace{5mm}

The idea of
finite dimensional approximation is
closely related to the notion of
``low energy effective theory" in physics.
Actually the approximation should be regarded
just as a part of the vast notion which
we can deal with rigorously or mathematically.

Since Witten's realization of
Donaldson theory as a TQFT,
the formal relation between
mathematically regorous definition of invariants
and their formal path integral expressions
has suggested many things.
For instance,
the well-definedness of the Donaldson invariant
is based on the fact that
the formal dimension of the moduli space
increases when the instanton number goes up.
This fact seems equivalent to the other fact that the
pure Yang-Mills theory is asymptitotically free
(for N=2 SUSY theory) and its renormalized theory
does exists.

In the case of the finite dimensional
approximations of Seiberg-Witten theory,
the suspension maps  give
relations between many choices of approximations.
If we use some generalized cohomology theories
to detect our invariants,
the suspension maps induces the Thom isomorphisms,
or integrations along fibers.
If we compare this setting with physics,
the family of integrations look quite
similar to the renormalization group.
It seems the Thom classes which play the role of
vacua.
In this sense, one could say that
the family of finite approximations
is a topological version of the renormlization group.
This topological setting is very limitted.
It, however, has one advantage:
Usually the path integral expression is supposed to
take values in real or complex numbers.
On the other hand our invariants could take
values in torsions.

Let us conclude this survey by giving
three open problems.
\begin{enumerate}
\item
What is the correct formulation of
the geography of spin $4$-manifolds with $b_1>0$?
(If the intersection on $H^1$ is complicated enough,
then $\sign(X)$ would have stronger restriction.)
\item
When an oriented closed $3$-manifold is a
link of isolated algebraic singular point,
construct a canonical Galois group action on
some completion of
the Floer homotopy type of Kronheimer-Manolescu.
(This problem was suggested by
a hand-written
manuscript by D.~Johnson  in which
Casson-type invariants were defined.)
\item
The Seiberg-Witten map is quadratic.
Extract non-topological
information from this structure.
(Is it possible to approach the 11/8-conjecture
from this point of view?)
\end{enumerate}

\label{lastpage}

\end{document}